\documentclass[graybox]{svmult}

\usepackage{type1cm}        %
\usepackage{makeidx}         %
\usepackage{graphicx}        %
\usepackage{multicol}        %
\usepackage[bottom]{footmisc}%

\usepackage{newtxtext}       %
\usepackage[varvw]{newtxmath}       %

\makeindex             %

\usepackage{tikz, pgfplots}
\pgfplotsset{compat=1.16}
\usetikzlibrary{calc}
\usepackage{tikz-3dplot}
\usetikzlibrary{3d, arrows.meta}
\usetikzlibrary{external}
\usepackage{standalone}
\usepackage{ifthen}
\usepackage{xspace}
\usepackage{mathtools}

\usepackage[nameinlink]{cleveref}
\providecommand{\texorpdfstring}[2]{#1}
\usepackage[section]{algorithm}
\usepackage{algorithmic}

\newcommand{\R}{\mathbb{R}} %

\newcommand{\half}{1/2}

\newcommand{\pmat}[1]{\begin{pmatrix} #1 \end{pmatrix}}

\newcommand{\grad}{\nabla}
\renewcommand{\div}{\nabla \cdot}
\newcommand{\dintx}{\, \mathrm{d} x}

\newcommand{\dints}{\, \mathrm{d} \sigma}

\newcommand{\dt}{\partial_t}
\newcommand{\dn}{\partial_n}
\newcommand{\interior}{\operatorname{int}}
\newcommand{\dom}{\operatorname{dom}}

\newcommand\restr[2]{{          %
\left.\kern-\nulldelimiterspace %
#1                              %
\vphantom{\big|}                %
\right|_{#2}                    %
}}

\newcommand{\norm}[1]{\left\lVert #1 \right\rVert}

\DeclarePairedDelimiterX\setc[2]{\{}{\}}{\,#1 \;\delimsize\vert\; #2\,}
\DeclarePairedDelimiter\monosetc{\{}{\}}

\newcommand{\Th}[1][]{
	\ifthenelse{ \equal{#1}{} }
	{\mathcal{T}_h}
	{\mathcal{T}_h(#1)}
}

\newcommand{\Fh}[1][]{
	\ifthenelse{ \equal{#1}{} }
	{\mathcal{F}_h}
	{\mathcal{F}_{h}(#1)}
}
\newcommand{\Fhint}[1][]{
	\ifthenelse{ \equal{#1}{} }
	{\mathcal{F}_h^\mathrm{int}}
	{\mathcal{F}_{h}^\mathrm{int}(#1)}
}
\newcommand{\Fhbnd}[1][]{
	\ifthenelse{ \equal{#1}{} }
	{\mathcal{F}_h^\mathrm{bnd}}
	{\mathcal{F}_{h}^\mathrm{bnd}(#1)}
}

\newcommand{\PkK}{\mathbb{P}_k(K)}

\newcommand{\PbkTh}[1][]{\mathcal{P}^b_k(\Th[#1])}

\newcommand{\ip}[2][]{
	\ifthenelse{ \equal{#1}{} }
	{\bigl( #2 \bigr)}
	{\bigl( #2 \bigr)_{#1}}
}
\newcommand{\abilSymb}{a} %
\newcommand{\abil}[2][]{
	\ifthenelse{ \equal{#1}{} }
	{\abilSymb\bigl( #2 \bigr)}
	{\abilSymb_{#1}\bigl( #2 \bigr)}
}

\newcommand{\Lh}[1][]{
	\ifthenelse{ \equal{#1}{} }
	{\mathcal{L}_h}
	{\mathcal{L}_{h,#1}}
} 
\newcommand{\LTproj}[1][]{
	\ifthenelse{ \equal{#1}{} }
	{\Pi_{\fulldomain}}
	{\Pi_{#1}}
} 

\newcommand{\nF}{n_F}
\newcommand{\hF}{h_{F}}

\newcommand{\av}[1]{\{ \! \! \{ #1 \} \! \! \}}
\newcommand{\weight}{\omega}
\newcommand{\avw}[1]{\av{#1}^{\weight}}
\newcommand{\jp}[1]{\llbracket #1 \rrbracket}
\newcommand{\njp}[1]{\cdot \nF \jp{#1}}

\newcommand{\dGpen}{\eta}

\newcommand{\timestep}{time step\xspace}
\newcommand{\timesteps}{time steps\xspace}
\newcommand{\timestepsize}{time-step size\xspace}

\renewcommand{\ts}[1][]{%
	\ifthenelse{ \equal{#1}{} }
	{\ensuremath{\tau}}
	{\ensuremath{\frac{\tau}{#1}}}
}
\newcommand{\tss}[1][]{%
	\ifthenelse{ \equal{#1}{} }
	{\ensuremath{\tau^2}}
	{\ensuremath{\frac{\tau^2}{#1}}}
}

\newcommand{\genD}{\widehat{\fulldomain}}

\newcommand{\genF}{\widehat{\Gamma}}
\newcommand{\genFh}{\Fh[\genF]}
\newcommand{\genDB}{\widehat{\Gamma}_D}
\newcommand{\genNB}{\widehat{\Gamma}_N}

\newcommand{\Patch}[2]{\mathcal{N}_{#1}(#2)}

\newcommand{\ov}{\delta}
\newcommand{\ovp}{\ell}
\newcommand{\fulldomain}{\Omega}
\newcommand{\dfulldomain}{\partial \Omega}
\newcommand{\fulldomainc}{\overline{\fulldomain}}
\newcommand{\SD}[1]{\fulldomain_{#1}}

\newcommand{\SDc}[1]{\overline{\fulldomain}_{#1}}

\newcommand{\SDov}[1]{\fulldomain_{#1}^{\ovp}}
\newcommand{\dSDov}[1]{\partial \SDov{#1}}
\newcommand{\SDovc}[1]{\overline{\fulldomain}_{#1}^{\ovp}}

\newcommand{\numberSD}{\mathcal{I}}
\newcommand{\SDovint}[1]{\Gamma_{#1}^{\ovp}}
\newcommand{\SDint}[1]{\Gamma_{#1}}
\newcommand{\PSD}[1]{\SDov{\SDint{},#1}}
\newcommand{\PSDc}[1]{\SDovc{\SDint{},#1}}

\newcommand{\sumSD}{\sum_{i=1}^{\numberSD}}

\newcommand{\solutionVector}{x}  
\newcommand{\uComponent}{u}
\newcommand{\vComponent}{v} 

\newcommand{\solDS}[1]{\solutionVector_{\normalfont\text{DS}}^{#1}} 
\newcommand{\uDS}[1]{\uComponent_{\normalfont\text{DS}}^{#1}}
\newcommand{\vDS}[1]{\vComponent_{\normalfont\text{DS}}^{#1}}

\newcommand{\solDSloc}[2]{\solutionVector^{#2}_{#1}} 
\newcommand{\uDSloc}[2]{\uComponent_{#1}^{#2}}
\newcommand{\vDSloc}[2]{\vComponent_{#1}^{#2}}

\newcommand{\uPred}[1]{\uComponent_{\star}^{#1}}
\newcommand{\vPred}[1]{\vComponent_{\star}^{#1}}

\begin{document}

\title*{A non-iterative domain decomposition time integrator combined with discontinuous Galerkin space discretizations for acoustic wave equations}

\titlerunning{Discontinuous Galerkin domain splitting for acoustic wave equations}

\author{Tim Buchholz\orcidID{0009-0000-0389-0983} and\\ Marlis Hochbruck\orcidID{0000-0002-5968-0480}}

\institute{Tim Buchholz \at Institute for Applied and Numerical Mathematics,
    Karlsruhe Institute of Technology, Englerstr. 2, 76131 Karls\-ru\-he, Germany, \email{tim.buchholz@kit.edu}
\and Marlis Hochbruck \at Institute for Applied and Numerical Mathematics,
    Karlsruhe Institute of Technology, Englerstr. 2, 76131 Karls\-ru\-he, Germany \email{marlis.hochbruck@kit.edu}}

\maketitle

\newcommand{\abstracttext}{
    We propose a novel non-iterative domain decomposition time integrator for acoustic wave equations using a discontinuous Galerkin discretization in space. It is based on a local Crank-Nicolson approximation combined with a suitable local prediction step in time. In contrast to earlier work using linear continuous finite elements with mass lumping, the proposed approach enables higher-order approximations and also heterogeneous material parameters in a natural way.
}

\abstract*{\abstracttext}
\abstract{\abstracttext}

\setlength{\belowdisplayskip}{6pt} \setlength{\belowdisplayshortskip}{6pt}
\setlength{\abovedisplayskip}{6pt} \setlength{\abovedisplayshortskip}{6pt}
\setlength{\intextsep}{10pt}

\section{Introduction}
\label{sec:Introduction}
    We construct a novel non-iterative domain decomposition time integrator for acoustic wave equations which uses a discontinuous Galerkin (DG) space discretization. Employing a DG discretization offers two key advantages. First, it easily allows us to use higher-order polynomials on the mesh elements and second, it works very well for 
    spatially varying material parameters, e.g., piecewise constant material coefficients modeling composite materials.
    The construction is inspired by the work of \cite{BluLR92,DawD92} for parabolic problems and
    \cite{BucH25}, where we proposed and analyzed this method for linear acoustic wave equations using a space discretization based on linear finite elements combined with mass lumping. 

    The linear acoustic wave equation is posed on an open, bounded, and polyhedral domain $\fulldomain \subset \mathbb{R}^d$ with a non-empty Dirichlet boundary $\Gamma_D \subseteq \dfulldomain$ and Neumann boundary $\Gamma_N = \dfulldomain \setminus \Gamma_D$. 
    The material coefficient $\kappa \in L^{\infty}(\fulldomain)$ satisfies $\alpha < \kappa(x) < \beta$ almost everywhere for some constants $\alpha, \beta > 0$, and may in particular be piecewise constant.
    Let $L u = \div( \kappa \grad u)$ be the differential operator applied to a function $u$ in the domain $D(L) = H^1(\fulldomain) \cap \setc{u \in L^2(\fulldomain)}{Lu \in L^2(\fulldomain)}$.
    Given initial data $u^0 \in D(L)$ and $v^0 \in H^1(\fulldomain)$, the linear acoustic wave equation is given by
    \begin{subequations}
        \begin{alignat}{6}
                \dt u &= v , \quad 
                    &\dt v &= Lu + f, 
                    \quad &&\text{in } \fulldomain \times (0,T] \\
                u(x,0) &= u^0(x), \quad 
                    &v(x,0) &= v^0(x), 
                    &&\text{in } \fulldomain \\
                && u &= g,
                    &&\text{on } \Gamma_D \times (0,T]
                    \\
                &&\dn u &= 0,
                    &&\text{on } \Gamma_N \times (0,T]
        \end{alignat}
    \end{subequations}
    where $T>0$ denotes the final time.
    For the inhomogeneity $f$ we assume $f \in C([0,T], H^1(\fulldomain))$. However, the precise necessary conditions on $g$ to get well-posedness of the problem are quite delicate, and lie outside the scope of this work. We just demand $g \in C^2([0,T], C^2(\Gamma_D))$, which is sufficient to lift the problem into one with homogeneous mixed boundary conditions. Moreover, 
    we assume compatibility of $g$ at the transition between $\Gamma_D$ and $\Gamma_N$, as well as $\restr{u^0}{\Gamma_D} = \restr{g}{t=0}$. 

    The remainder of the paper is structured as follows. In \Cref{sec:Preliminaries}, we review relevant preliminaries, including the space discretization, global time integrators, and cell extensions within a mesh. \Cref{sec:domainsplitting} presents the construction and details of the proposed domain splitting method.
    In \Cref{sec:implementation}, we highlight some key aspects to be considered for an efficient implementation. 
    Finally, \Cref{sec:numericalexp} presents numerical experiments that demonstrate the method's performance.
\section{Preliminaries}
\label{sec:Preliminaries}
    \subsection{Discretization in space}
        For the spatial discretization, we consider a shape- and contact-regular, matching simplicial mesh $ \Th = \Th[\fulldomain] $ of the domain $ \fulldomain $, cf.\ \cite[Definition 8.11]{ErnG21I}.
        We denote by $\Fh = \Fh[\fulldomain]$ the set of all mesh faces, decomposed into the set of boundary faces $\Fhbnd = \Fhbnd[\fulldomain]$ and the set of interior faces $\Fhint = \Fhint[\fulldomain]$.
        Further, we assume that $\dfulldomain$ is the distinct union of the Dirichlet boundary $\Gamma_{D}$ and the Neumann boundary $\Gamma_{N}$ and that all faces are contained either in $\Gamma_D$ or in $\Gamma_N$.
        The wave propagation speed $\kappa$ is piecewise constant, and the mesh $\Th$ is matched to $\kappa$, i.e., $\restr{\kappa}{K}$ is constant for every element $K \in \Th$.  
        For each element $K \in \Th$, let $h_K$ denote its diameter and $h$ the maximal element diameter in the mesh.

        For any subset $\widehat{\Th} \subset \Th$ of the mesh, we define the corresponding spatial domain
        \begin{equation}
            \label{eq:DomTh}
            \genD =	\dom(\widehat{\Th}) \coloneqq \interior  \bigcup_{K \in \widehat{\Th}} \overline{K} \subset \Omega.
        \end{equation}
        If a domain $\genD$ is defined in this way as a union of cells, then the set of mesh elements belonging to $\genD$ is denoted by 
        $\Th[\genD] = \{ K \in \Th \, \mid \, K \subset \genD \}.$
        We further denote by $\Fh[\genD]$ the set of all faces in $\Th[\genD]$, which we split into the interior faces $\Fhint[\genD]$ and boundary faces $\Fhbnd[\genD]$.
        The Dirichlet and Neumann parts of the boundary of $\genD$ are denoted by $\genDB$ and $\genNB$, respectively.  
        Analogously, for a generic interface $\genF$ we denote its associated set of faces by $\genFh \subset \Fh$. As defined in \cite[Definition~1.18]{DiPE12} we denote the face normals of a given face $F\in \Fh$ by $n_F$.

        Next, we introduce the discrete function spaces.  
        For each element $K \in \Th$, let $\PkK$ denote the set of all polynomials in $d$ variables of total degree at most $k$.  
        The associated broken polynomial space on a domain $\genD$ is then defined as
        \begin{equation*}
                \PbkTh[\genD] \coloneqq \setc{\psi_h \in L^2(\genD; \R)}{\restr{\psi_h}{K} \in \PkK, \; \forall K \in \Th[\genD]}.
        \end{equation*}
        This broken space serves as the approximation space in the discontinuous Galerkin method, 
        allowing discontinuities of the functions across element interfaces.
        We denote the standard $L^2$ inner product on $\genD$ by $ \ip[\genD]{\cdot,\cdot} = \ip[L^2(\genD)]{\cdot,\cdot}$.
        The corresponding $L^2$ projection $\Pi_{\genD}: L^2(\genD) \to \PbkTh[\genD]$ is defined for any $\phi \in L^2(\genD)$ by
        \begin{equation}
            \label{eq:defL2project}
            \ip[\genD]{\LTproj[\genD] \phi - \phi, \psi_h} = 0 , \quad \forall \psi_h \in \PbkTh[\genD].
        \end{equation}
        We follow the symmetric weighted interior penalty (SWIP) approach as described in 
        \cite[Section~4.5.2.3]{DiPE12}, where jumps are denoted by $\jp{\cdot}$ and weighted averages by 
        $\avw{\cdot}$ with weights defined as in \cite[Definition~4.46]{DiPE12}. 
        In this context, we also recall the 
        definition of the local length scale $h_F$, cf.\ \cite[Definition~4.5]{DiPE12}, the 
        $\kappa$-dependent penalty parameter $\gamma_{\kappa,F}$ introduced in \cite[below eq.~(4.64)]{DiPE12}, 
        and the penalty parameter $\eta > 0$.
        Denoting the broken gradient from \cite[Definition~1.21]{DiPE12} by $\nabla_h$ the SWIP bilinear form on $\genD$ is given by
        \begin{align*} 
            \abil[\genD]{\phi_h,\psi_h} 
                 =&  \sum_{K \in \Th[\genD]}\int_{K} \kappa \grad_h \phi_h \grad_h \psi_h \dintx + \sum_{F \in \Fh[\genD]}  \dGpen\frac{ \gamma_{\kappa,F}}{\hF}  \int_F \jp{\phi_h} \jp{\psi_h} \dints \\
                    &  - \sum_{F \in \Fh[\genD]} \int_F \avw{\kappa \grad_h \phi_h} \njp{\psi_h} + \avw{\kappa \grad_h \psi_h} \njp{\phi_h} \dints ,
        \end{align*}
        for $\phi_h, \psi_h \in \PbkTh[\genD]$. We write $\norm{\cdot}_{a,\genD}^2$ for the associated SWIP norm, cf. \cite[eq.~(4.69)]{DiPE12}.
        We now introduce the linear operator 
        $\Lh[\genD]: \PbkTh[\genD] \to \PbkTh[\genD]$ 
        associated with the SWIP bilinear form. It is defined by 
        \begin{equation}
            \label{def::Lh}
            \ip[\genD]{\Lh[\genD] \phi_h, \psi_h} = \abil[\genD]{\phi_h,\psi_h} , \quad \forall \phi_h,\psi_h \in \PbkTh[\genD] .
        \end{equation}
        For the case $\genD = \fulldomain$, we write simply $\Lh = \Lh[\fulldomain]$.
        For the treatment of inhomogeneous Dirichlet boundary conditions we refer to 
        \cite[Section~4.2.2]{DiPE12}.
        Specifically, the boundary data $g$ on a set of faces $\genF$ is weakly enforced 
        by an additional term on the right-hand side of the discrete problem:
        We thus define
        \begin{equation}
            \label{eq:DefBcTerm}
            \ip[\genD]{\mathcal{G}_{\genF}(g), \psi_h}
            =  \sum_{F \in \genFh}
                 \dGpen\frac{ \gamma_{\kappa,F}}{\hF} \int_{F} g \psi_h \dints
                 -  \int_{F} g \kappa \grad_h \psi_h \cdot n \dints,
        \end{equation}
        for all $\psi_h \in \PbkTh[\genD] $.
        We abbreviate  $\mathcal{G}_{\genF}(g^n) \coloneqq \mathcal{G}_{\genF}(g(t_n))$ for a given $t_n \in [0,T]$.

        The initial values and the right-hand side are approximated using the $L^2$ projection onto the broken polynomial space. 
        Specifically, on a domain $\genD$ we set
        \begin{equation}  \label{eq:projected-data}
            u_{h,\genD}^0 = \LTproj[\genD] u^0, 
            \qquad 
            v_{h,\genD}^0 = \LTproj[\genD] v^0, 
            \qquad 
            f_{h,\genD}(t) = \LTproj[\genD] f(t).
        \end{equation}
        Whenever the domain is clear from the context, we omit the subscript $\genD$, and we also abbreviate 
        $f_h^n = f_{h,\genD}(t_n)$ for a given $t_n \in [0,T]$.
        With these definitions, the semi-discrete problem in 
        $\PbkTh[\genD] \times \PbkTh[\genD]$ reads
        \begin{equation}
            \label{eq:semiDiscreteProblem}
            \partial_t \pmat{u_h \\ v_h} = \pmat{0 & I \\ - \Lh[\genD] & 0 } \pmat{u_h \\ v_h} + \pmat{0 \\ f_{h,\genD}(t)} + \pmat{0 \\ \mathcal{G}_{\genDB}(g(t))}.
        \end{equation}

    \subsection{Time integrators on \texorpdfstring{$\genD$}{a generic subdomain}}
        We consider a generic subdomain $\genD$ and a given \timestepsize $\ts > 0$ with
        \[
            t_n = n \ts, \qquad n = 0, \dots, N_T, \quad T = N_T \ts .
        \]
        Starting from the semi-discrete problem \eqref{eq:semiDiscreteProblem}, we present two standard second-order accurate time integration methods to obtain full-discretizations on $\genD\times[0,T]$.
        As described in \Cref{sec:Preliminaries} the boundary data $g(t_n)$ is incorporated weakly through the
        operator  $\mathcal{G}_{\genDB}$ defined in \eqref{eq:DefBcTerm}.  
        The \textbf{Crank-Nicolson} discretization of \eqref{eq:semiDiscreteProblem} is given by
        \begin{subequations}  
            \label{eq:CN}
            \begin{align}
                \label{eq:CN-u}
                    u_h^{n+1} &= u_h^{n} + \ts[2] 
                    \bigl( 	v_h^{n+1}  + 	v_h^{n}  \bigr),\\
                \label{eq:CN-v}
                    v_h^{n+1} &= v_h^{n} - \ts[2] \Lh[\genD] \bigl( 	u_h^{n+1}  + 	u_h^{n}  \bigr)
                    + \ts[2]  \bigl( 	 f_h^{n+1} + 	f_h^{n} \bigr) \\
                    \nonumber 
                    & \quad + \ts[2]  \bigl( 	 \mathcal{G}_{\genDB}(g^{n+1}) + 	\mathcal{G}_{\genDB}(g^{n}) \bigr).
            \end{align}
        \end{subequations}
        Moreover, we also consider the \textbf{leapfrog} method, which reads
        \begin{subequations}  
            \label{eq:LF}
            \begin{align}
                \label{eq:LF-vhalf}
                    v_h^{n+\half} &= v_h^n - \ts[2] \Lh[\genD]  u_h^n + \ts[2]f_h^n + \ts[2]\mathcal{G}_{\genDB}(g^{n})\\
                \label{eq:LF-u}
                    u_h^{n+1} &= u_h^n + \ts v_h^{n+\half} \\
                \label{eq:LF-v}
                    v_h^{n+1} &= v_h^{n+\half} - \ts[2] \Lh[\genD] u_h^{n+1} + \ts[2] f_h^{n+1} + \ts[2] \mathcal{G}_{\genDB}(g^{n+1}).
            \end{align}
        \end{subequations}
        In the discontinuous Galerkin setting, the leapfrog method has the advantage that the resulting mass matrix %
        is block-diagonal. There, each block in the mass matrix corresponds to the degrees of freedom in a single cell $K \in \Th$.

        \subsection{Cell extensions}

        \begin{figure}
                \sidecaption[t]
                \input{overlap.tex}
                \caption{Extension $\mathcal{N}_2(\genD)$ by $\ell = 2$ layers of a subdomain $\genD \subset \fulldomain$ colored in dark blue. 
                }
                \label{Fig:CellExtension}
        \end{figure}

        To define overlapping subdomains and cell neighborhoods, we introduce the concept of cell extensions.  
        For a given subdomain $\genD \subset \fulldomain$ and a fixed integer $\ell\geq 1$, we define its extension by $\ell$ layers of neighboring cells recursively as
        \begin{equation}
            \label{eq:defPatches}
            \begin{aligned}
                \Patch{0}{\genD} &\coloneqq  \Th[\genD], \\
                \Patch{j}{\genD} &\coloneqq \setc{K \in \Th[\fulldomain]}{\exists K_{\star} \in \Patch{j -1}{\genD} \; : \; \overline{K} \cap \overline{K}_{\star} \neq \emptyset},
            \end{aligned}
        \end{equation}
        for $j=1,\ldots,\ell$, see also \Cref{Fig:CellExtension}.
        Similarly, for a generic interface $\genF$, we define the interface cell extension as
        \begin{equation}
            \label{eq:ExtensionInterface}
            \begin{aligned}
                \Patch{}{\genF} &\coloneqq \setc{K \in \Th[\genD]}{\exists F_{\star} \in \genFh \; : \; \overline{K} \cap F_{\star} \neq \emptyset},
            \end{aligned}
          \end{equation}
          which contains all elements which share a face or a corner with $\genF$.

\section{Domain splitting method}
    \label{sec:domainsplitting}
    To construct our new method, we first decompose the spatial domain $\fulldomain$ into $\numberSD$ distinct, non-overlapping subdomains $\SD{i}$, i.e.
    \begin{equation}
        \label{eq:NonOvDecomp}
        \fulldomainc = \bigcup\limits_{i=1}^{\numberSD} \SDc{i},
        \qquad
        \Omega_i \cap \Omega_j = \emptyset \quad \text{for } i \neq j
        \;.
    \end{equation}
    Based on the definitions \eqref{eq:DomTh} and \eqref{eq:defPatches}, we then introduce overlapping subdomains $\SDov{i}$ by extending each $\SD{i}$ with $\ell$ layers of elements 
    \begin{equation}
        \label{eq:OvDecomp}
        \SDov{i} \coloneqq \dom \Patch{\ell}{\SD{i}}, \quad i = 1, \dots \numberSD.
    \end{equation}
    Similarly, using \eqref{eq:ExtensionInterface}, we define domains based on the cells around the interface  $\SDovint{i} = \dSDov{i}\cap \fulldomain$
    \begin{equation}
                \label{eq:PSD}
                \PSD{i} \coloneqq \dom \Patch{}{\SDovint{i}}, \quad i = 1, \dots \numberSD,
    \end{equation}
    which we will refer to as prediction domain for the %
    interface $\SDovint{i}$ of $\SDov{i}$. An overview of this subdomain notation is given in \Cref{Fig:MeshesImplementation}.
    \begin{figure}
        \centering
        \input{overlapping_domain.tex}
        \input{prediction_domain.tex}
        \caption{Overlapping subdomain $\SDov{i}$ (left, dark and light red area) and prediction domain $\PSD{i}$ (right, yellow area).
        The interface $\SDovint{i}$ is shown in dark red in both pictures.
        }
        \label{Fig:MeshesImplementation}
    \end{figure}
    Let $\ov$ be the minimal physical width of the overlap $\SDov{i} \setminus \SDc{i}$, which satisfies $\ov \sim \ell h$, given the underlying mesh $\Th $ is shape- and contact-regular. We denote by $\Gamma^{\ell}_{i,D}$ the portion of the boundary of $\SDov{i}$ that coincides with the Dirichlet boundary of $\fulldomain$, i.e., 
    \[
        \Gamma^{\ell}_{i,D} = \Gamma_D \cap \SDovc{i}.
    \]

    Next, we describe the construction of the domain-splitting approximations
    \begin{equation}
        \solDS{n} = \pmat{\uDS{n} \\ \vDS{n}}, \qquad \text{for } n=0,\dots,N_T.
    \end{equation}
    For the initial values we use $L^2$ projections
    \begin{equation}
        \label{eq:InitialValuesDS}
        \uDS{0} = u_h^0 = \LTproj[\fulldomain] u^0 ,\qquad  \vDS{0} = v_h^0 = \LTproj[\fulldomain] v^0 .
    \end{equation}
    Note, that we can also replace this by local $L^2$ projections on the subdomains
    \begin{equation*}
        \restr{\uDS{0}}{\SDovc{i}} 
        = \restr{\bigl(\LTproj[\fulldomain] u^0\bigr)}{\SDovc{i}} 
        = \LTproj[\SDov{i}] u^0 \; ,
    \end{equation*}
    since we are using a discontinuous Galerkin discretization.
    Given an approximation $\solDS{n}$ at time $t_n$, the method advances to the next time step as described below:

    First, we loop over $i=1, \dots, \numberSD$ and perform for each subdomain:

    \begin{trailer}{Prediction on strip $\PSD{i}$ of the interface $\SDovint{i}$}%
        A leapfrog step is carried out on $\PSD{i}$ to compute an approximation $\uPred{n+1}$:
    \begin{subequations}  
        \label{eq:LFpred}
        \begin{align}
            \label{eq:LFpred-vhalf}
                \vPred{n + \half} &= \restr{\vDS{n}}{\PSDc{i}} - \ts[2] \Lh[\PSD{i}]  \restr{\uDS{n}}{\PSDc{i}} + \ts[2]\restr{f_h^n}{\PSDc{i}} + \ts[2]\mathcal{G}_{\Gamma_D \cap \PSDc{i}}(g^{n})\\
            \label{eq:LFpred-u}
                \uPred{n+1} &= \restr{\uDS{n}}{\PSDc{i}} + \ts \vPred{n+\half}
        \end{align}
    \end{subequations}
    Using $\uPred{n+1}$ we define a boundary term $\mathcal{G}_{\SDovint{i}}(\avw{\uPred{n+1}})$ on the interface $\SDovint{i}$.
    \end{trailer}
    \newpage
    \begin{trailer}{Local calculation on $\SDov{i}$}%
        On $\SDov{i}$, we perform a Crank-Nicolson step with weakly imposed interface and boundary data.
        This results in a boundary term of the form
    \begin{subequations}  
        \label{eq:locCN}
        \begin{equation}
        \label{eq:locCN-bc}
        \mathcal{G}_{\SDov{i}}^{n+1,n} \coloneqq \ts[2] \Bigl( 
            \mathcal{G}_{\SDovint{i}}(\avw{u^{n+1}_{\star}}) 
            + \mathcal{G}_{\Gamma^{\ell}_{i,D}}(g^{n+1}) 
            + \mathcal{G}_{\SDovint{i}}(\avw{\uDS{n}}) 
            + \mathcal{G}_{\Gamma^{\ell}_{i,D}}(g^n) 
            \Bigr).
        \end{equation}
        The resulting subdomain approximations $\solDSloc{i}{n+1}$ on $\SDovc{i}$ are then given by
        \begin{align}
            \label{eq:locCN-u}
                \uDSloc{i}{n+ 1} &= \restr{\uDS{n}}{\SDovc{i}} + \ts[2] 
                \bigl( 	\vDSloc{i}{n+1}  + 	\restr{\vDS{n}}{\SDovc{i}} \bigr),\\
            \label{eq:locCN-v}
                \vDSloc{i}{n + 1} &= \restr{\vDS{n}}{\SDovc{i}} - \ts[2] \Lh[\SDov{i}] \bigl( 	\uDSloc{i}{n+1}  + 	\restr{\uDS{n}}{\SDovc{i}}  \bigr)
                + \ts[2]  \restr{\bigl( 	 f_h^{n+1} + 	f_h^{n} \bigr)}{\SDovc{i}}
                + \mathcal{G}_{\SDov{i}}^{n+1,n}.
        \end{align}
    \end{subequations}
    \end{trailer}
    Note, that \eqref{eq:LFpred} and \eqref{eq:locCN} do not couple between subdomains and can thus both be performed in parallel across all subdomains.  
    Then, after completion of the loop over all subdomains we exchange data among neighboring subdomains:

    \begin{trailer}{Exchange approximations in overlap regions}%
            If two subdomains $\SD{i}$, $\SD{j}$ are adjacent, i.e., 
            $S_{ij}= \SDov{i} \cap \SDov{j} \neq \emptyset$, then the local approximations are exchanged across the overlap. For the efficiency, it is important that this does not require global communication.
            Specifically, we replace
            \[
                \restr{\solDSloc{i}{n+1}}{S_{ij}\cap {\SDc{j}}}  \quad \text{by} \quad  \restr{\solDSloc{j}{n+1}}{S_{ij}\cap {\SDc{j}}}
                \qquad \text{and} \qquad 
                \restr{\solDSloc{j}{n+1}}{S_{ij}\cap {\SDc{i}}} \quad \text{by} \quad \restr{\solDSloc{i}{n+1}}{S_{ij}\cap {\SDc{i}}} \; .
            \]
            We store these updated values in
            \begin{equation}
                \label{eq:localUpdates}
                \restr{\solDS{n+1}}{\SDovc{i}} \quad \text{and} \quad \restr{\solDS{n+1}}{\PSDc{i}}
            \end{equation}
            locally on each subdomain for $i = 1,\dots, \numberSD$ for the next \timestep.
    \end{trailer}

    These local updates are equivalent to assembling the global approximation via
    \begin{equation}
        \label{eq:globalApprox}
        \solDS{n+1} 
        \coloneqq \sumSD \restr{\solDSloc{i}{n+1}}{\SDc{i}},
    \end{equation}
    where we sum up the local approximations restricted to $\SDc{i}$.
    Note however, that it is not necessary to act globally here.
    Moreover, unlike in the mass-lumped finite element setting (cf.\ \cite{BucH25}), no averaging at the non-overlapping interfaces is required. We summarize the method in \Cref{alg:DSDG}.
    \begin{algorithm}
    \caption{Domain splitting for DG discretizations (one \timestep)}
    \label{alg:DSDG}
    \begin{algorithmic}
    \STATE  $\{\text{given } \restr{\solDS{n}}{\SDovc{i}} \text{ and } \restr{\solDS{n}}{\PSDc{i}}\}$ 
    \FOR{$i=1:\numberSD$}
    \STATE $\{ \text{on each subdomain}\}$
    \STATE calculate prediction $u^{n+1}_{\star}$ by leapfrog step \eqref{eq:LFpred} on $\PSD{i}$
    \STATE calculate $\solDSloc{i}{n+1}$ by Crank-Nicolson step \eqref{eq:locCN} on $\SDov{i}$ with boundary term $\mathcal{G}_{\SDov{i}}^{n+1,n}$ \eqref{eq:locCN-bc}
    \ENDFOR
    \STATE update overlap regions with values from  $\restr{\solDSloc{i}{n+1}}{\SDc{i}}$, $i = 1, \dots, \numberSD$
    \STATE $\to $ prepare $\restr{\solDS{n+1}}{\SDovc{i}}$ and $\restr{\solDS{n+1}}{\PSDc{i}}$ for the next step
    \STATE $\{  \text{only if desired build a global approximation } \solDS{n+1} \text{ on } \fulldomain \text{ by \eqref{eq:globalApprox}}\}$
    \end{algorithmic}
    \end{algorithm}
\section{Implementation}
    \label{sec:implementation}
    The method is implemented within the FEniCSx framework \cite{BarDDetAl23}, using its PETSc interface via petsc4py \cite{DalPetAl11}. 
    The global mesh is generated with Gmsh \cite{GeuRe09}, which is subsequently partitioned into non-overlapping subdomains using PT-Scotch \cite{CheP08}, cf. \Cref{Fig:ScotchPartitioning}.
    For each subdomain, we construct the corresponding local meshes of the overlapping subdomain $\SDov{i}$ and the leapfrog prediction domain $\PSD{i}$, as illustrated in \Cref{Fig:MeshesImplementation}. The implementation is consistently based on these local meshes $\Th[\SDov{i}]$ and $\Patch{}{\SDovint{i}}$, while direct access to the global mesh $\Th[\fulldomain]$ is avoided whenever possible.

    Communication between subdomains is performed using point-to-point MPI routines. Since a subdomain can only communicate with one neighbor at a time, the exchanges must be organized into rounds in order to avoid deadlocks. To structure these rounds, we construct a weighted undirected graph whose nodes correspond to subdomains and whose edges represent the messages (value exchanges) between neighboring subdomains; the weights encode the number of values to be transmitted, see \Cref{Fig:scheduling}.

    \begin{figure}
        \sidecaption[t]
        \input{scheduling.tex}
        \caption{Example for a weighted undirected communication graph for $\numberSD = 6$ subdomains. Weights correspond to the number of values, which need to get exchanged between two subdomains.}
        \label{Fig:scheduling}
	\end{figure}
 
    The sequence of communication rounds is determined by a greedy algorithm with two objectives: minimizing the total number of rounds and balancing the message sizes within each round. To this end, the edges are first sorted by message size and then by the maximal node degree (i.e., the largest number of edges incident to either endpoint). Rounds are then built by selecting edges from the front of the sorted list, ensuring that no node appears more than once per round. This procedure yields a communication schedule that is both efficient and well-balanced.

    For the example graph in \Cref{Fig:scheduling}, the algorithm produces the following sequence of communication rounds:
    \begin{align*}
        &1: \hspace*{-1cm} &&(3,5,130), \hspace*{-1.75cm}  &&(2,4,90)   \\
        &2: \hspace*{-1cm} &&(1,3,120), \hspace*{-1.75cm}  &&(4,6,85)   \\
        &3: \hspace*{-1cm} &&(3,4,100), \hspace*{-1.75cm}  &&(5,6,110) ,  \;  (1,2,100) \\ 
        &4: \hspace*{-1cm} &&(4,5,20),  \hspace*{-1.75cm}  &&(2,3,15) \\
        &5: \hspace*{-1cm} &&(1,4,10),  \hspace*{-1.75cm} &&(3,6,10) 
    \end{align*}
    Here, a tuple $(i,\!j,\!M)$ denotes an exchange of $M$ values between subdomains $i$ and $j$.

    The exchange of values in the overlap regions requires particular care in the implementation. 
    To organize this process, we construct a dofmap that associates local DoFs with their corresponding global DoFs in $\Th[\fulldomain]$, see \Cref{Fig:dofmap}. 
    \begin{figure}
        \sidecaption
        \input{dofmap.tex}
        \caption{For each local mesh we build a 'dofmap', which maps the local DoFs to the global DoFs in $\Th[\fulldomain]$. Here the index denotes the local DoF, while the value stores the global DoF. Note, that the size of the map is only related to the local mesh.}
        \label{Fig:dofmap}
    \end{figure}
    This mapping allows us to identify exactly which DoFs in a subdomain require updated values from which neighbor. Moreover, it enables us to directly relate the DoFs of one subdomain to the DoFs of another subdomain, so that the communication can be set up locally and efficiently between neighboring subdomains.
    
    These preparatory steps --- global meshing, partitioning, dofmap setup, and scheduling --- are carried out once in a preprocessing phase. During the actual simulation, communication then proceeds directly between neighboring subdomains, without any further global operations on $\Th[\fulldomain]$, except when a global reconstruction as in \eqref{eq:globalApprox} is explicitly required.
\section{Numerical experiments}
    \label{sec:numericalexp}

    We simulate the propagation of a linear wave crossing a triangular inclusion, representing the cross-section of a prism, in two dimensions.
    The example features a spatially varying coefficient $\kappa(x)$ and mixed inhomogeneous boundary conditions.
    The computational domain is $\fulldomainc = [0,8]\times[0,4] \subset \R^2$, with a piecewise constant coefficient $\kappa$ taking values $\kappa_i = 1.0$ inside of the prism and $\kappa_o = 1.5$ outside, corresponding to a refractive index of $1.5$ which is typical for glass.
    An inhomogeneous Dirichlet condition is imposed on the inflow boundary $\Gamma_D = \restr{\dfulldomain}{x = 0}$ to generate an incoming wave of frequency $\omega = 0.0125$, described by
    \begin{equation*}
        g(x,y,t) = \sin  \bigl( \frac{t}{\omega} \bigr) W(y), 
    \end{equation*}
    where $W \in C^2(\R)$ satisfies $W(y)=1$ for $y \in [1,3]$ and $W(y)=0$ for $y \in [0,0.5] \cup [3.5,4.0]$. 
    Homogeneous Neumann conditions $\partial_n u = 0$ are applied on $\Gamma_N = \dfulldomain \setminus \Gamma_D$.

    The mesh $\Th[\fulldomain]$ is generated in Gmsh \cite{GeuRe09} with a local refinement (factor $2$) around the prism and partitioned into eight subdomains using PT-Scotch \cite{CheP08}, cf. \Cref{Fig:GmshPrism,Fig:ScotchPartitioning}.
    The simulation runs over $[0,T]$ with $T = 3.0$.

    \begin{figure*}
        \sidecaption
        \input{prism_h=0.25.tex}
        \caption{Example grid on $[0.8] \times [0,4]$ for the prism example with quite coarse $h$.  Red area $\kappa(x) = \kappa_i$, white area $\kappa = \kappa_o$. Dirichlet boundary $\Gamma_D$ (green) and Neumann boundary $\Gamma_N$ (blue). For simulations much finer $h$ were used.}
        \label{Fig:GmshPrism}
    \end{figure*}
    \begin{figure*}
        \sidecaption[t]
        \includegraphics[width=0.52\textwidth]{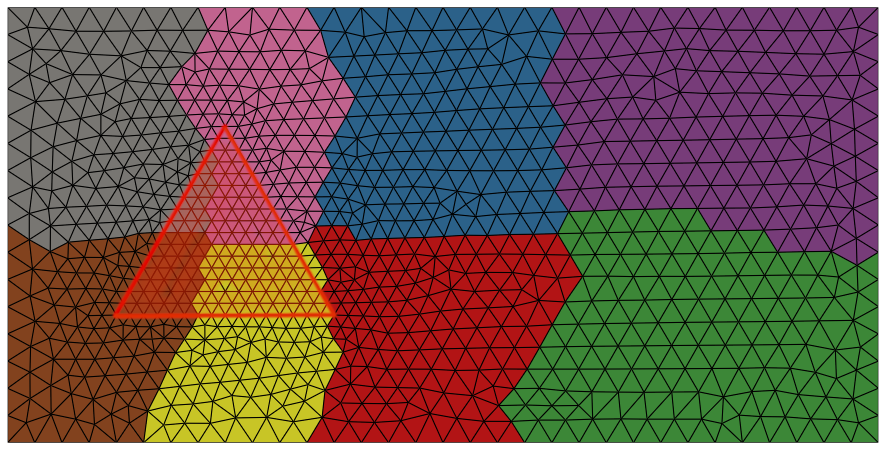}
        \caption{Mesh from \Cref{Fig:GmshPrism} partitioned into $8$ non-overlapping subdomains using the PT-Scotch partitioner \cite{CheP08}. The red area depicts the prism, where $\kappa(x) = \kappa_i$.
        }
        \label{Fig:ScotchPartitioning}
    \end{figure*}

    To evaluate the accuracy of the domain splitting method, we compare its results with a reference solution obtained from a leapfrog simulation on a refined mesh with $1{,}163{,}020$ cells. 
    A DG discretization with polynomial degree $p=2$ yields $6{,}978{,}120$ DoFs per component for the reference mesh.
    The domain splitting method is applied on a coarser mesh with $812{,}298$ cells ($4{,}873{,}788$ DoFs per component, $h_{\min}=0.00446$) using $8$ subdomains.
    We consider different overlap parameters $\ell \in \monosetc{2,4,8}$, denoted by $\text{DS}_2$, $\text{DS}_4$, and $\text{DS}_8$, respectively.
    For comparison, a global Crank-Nicolson (CN) simulation is performed on the same mesh. 
    The relative $L^2$-error of the $u$-component, $\norm{u_{\text{DS}} - u_{\text{lf}} }_{L^2(\fulldomain)} /\norm{u_{\text{lf}} }_{L^2(\fulldomain)}$ measured at $T=3.0$, is shown in \Cref{Fig:ConvergenceRef}. The results indicate that the CFL condition is linear dependent on $\ell$ analogously to the findings in \cite{BucH25}. However, this has to be rigorously analyzed in future work.

    \begin{figure*}
        \input{convergencePlot1.tex}
        \caption{Relative $L^2$-error of $u$ at final time $T=3.0$ for the domain splitting method ($\text{DS}_{\ell}$) with different overlap parameter $\ell$, compared to a reference leapfrog solution on a refined mesh. The Crank-Nicolson method ($\text{CN}$) is shown for comparison.}
        \label{Fig:ConvergenceRef}
    \end{figure*} 

    Next, we compare the domain splitting approximation directly to the global Crank-Nicolson (CN) solution.
    Both methods are run on the same mesh with $4{,}873{,}788$ DoFs and $h_{\min}=0.00446$.
    For the CN reference, we employ a Cholesky decomposition (chol), while the domain splitting systems are solved iteratively using a conjugate gradient (cg) method preconditioned with an incomplete Cholesky (icc) factorization.
    For this problem size, no speedup was observed when solving the global CN system iteratively.
    The relative difference between both solutions is measured in the combined $\norm{\cdot}_{a, \fulldomain} \times \norm{\cdot}_{L^2(\fulldomain)}$ norm; see \Cref{Fig:ConvergenceDiffCN}.

    \begin{figure*}
        \sidecaption[t]
        \input{convergencePlot2.tex}
        \caption{Difference of domain splitting approximation to Crank-Nicolson approximation at the end time $T=3.0$. We measure the first component in $\norm{\cdot}_{a,\fulldomain}$ and the second component in $\norm{\cdot}_{L^2(\fulldomain)}$.}
        \label{Fig:ConvergenceDiffCN}
    \end{figure*} 

\begin{figure*}
            \sidecaption[t]
            \includegraphics[width=0.3\textwidth]{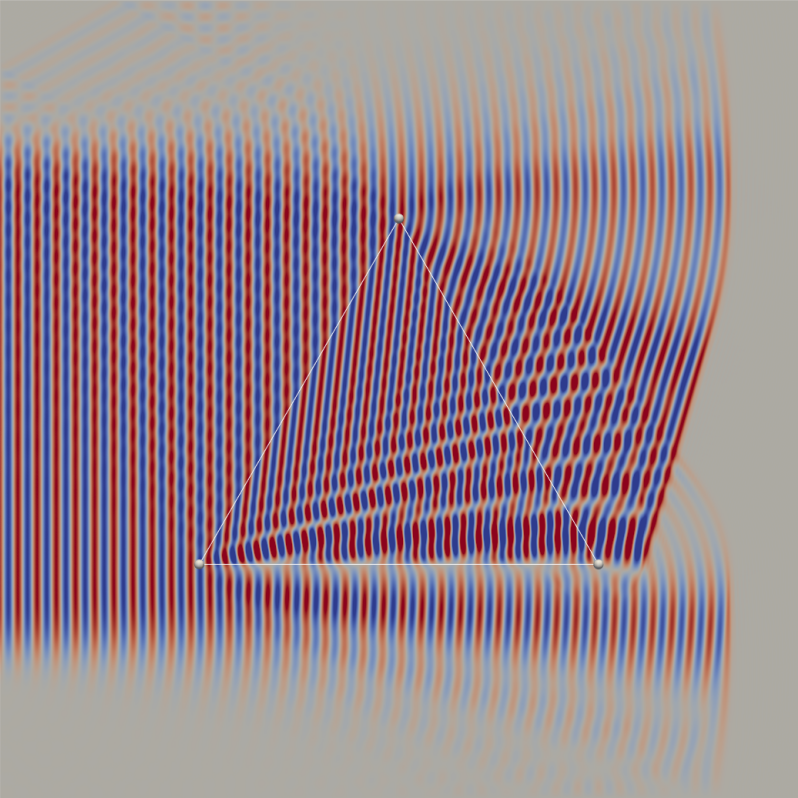} 
            \includegraphics[width=0.3\textwidth]{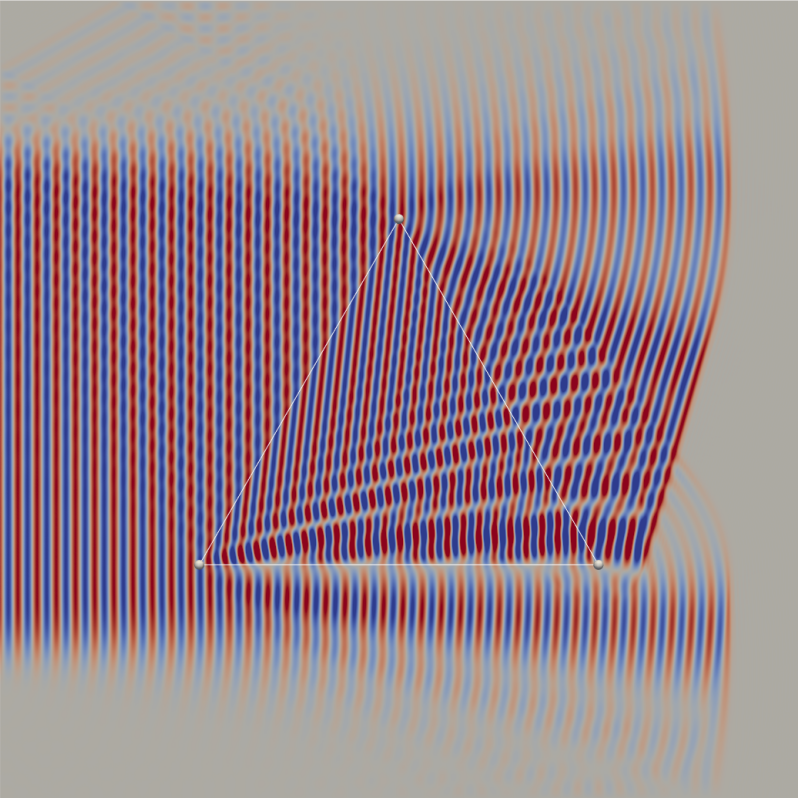}
            \caption{Snapshots of the domain splitting approximation $\text{DS}_4 $ (chol, left) and Crank-Nicolson approximation (chol, right) at the end time $T=3.0$ on the subarea $[0,4] \times [0,4]$.}
            \label{Fig:VizualizeApproximations}
    \end{figure*}

    Finally, \Cref{Table:comparisonDSCN} reports run times and solver performance for $\text{DS}_4$ and the global Crank-Nicolson (CN) method, parallelized via the MPI interfaces of FEniCSx and petsc4py.
    All simulations were carried out on the same workstation using $8$ cores, a time step $\tau = 0.001$, and the mesh from \Cref{Fig:ConvergenceDiffCN}.
    A visualization of both solutions at $T = 3.0$ is given in \Cref{Fig:VizualizeApproximations}.
    The results indicate that the proposed domain splitting method attains accuracy comparable to global time integration while providing structural benefits for parallelization.
    
    \setlength{\tabcolsep}{12pt}
     \begin{table}[H]
        \centering
        \begin{tabular}[t]{lccc}
        \hline
        &$\text{DS}_4$ (chol) & $\text{DS}_4$ (cg + icc)  & CN (chol)\\
        \hline
        rel. $L^2$ error in $u$ against ref &7.545e-2 &7.545e-2 & 7.547e-2\\
        rel. difference to CN (chol)&3.639e-5&3.639e-5 &--\\
        time for meshing  &94.9 s &94.9 s&70.8 s\\
        setup time for solvers &112.2 s &0.17 s & 47.9 s\\
        wall time per \timestep & 0.94 s & 1.21 s &1.52 s\\
        total simulation time & 3021 s & 3739 s  & 4669 s \\
        \hline
        \end{tabular}
        \caption{Detailed comparison between the domain splitting method with different solver configurations and the global (parallelized) Crank-Nicolson method for a simulation with $3000$ \timesteps.}
        \label{Table:comparisonDSCN}
    \end{table}%

    A more detailed theoretical study and large-scale experiments are planned as part of future work.
    The code corresponding to this section is made publicly available at 
    \begin{equation*}
        \text{\textcolor{blue!50!black}{\url{https://github.com/tim-buchholz/dsdg-acoustic-wave.git}}}.
    \end{equation*}

\begin{acknowledgement}
This work was funded by the Deutsche Forschungsgemeinschaft (DFG, German Research Foundation) — Project-ID 258734477 — CRC 1173.
The authors acknowledge support by the state of Baden-Württemberg through bwHPC. \\[0.5em]
\textbf{Competing Interests}
The authors have no conflicts of interest to declare that are relevant to the content of this chapter.
\end{acknowledgement}
\bibliographystyle{spmpsci}
\bibliography{references.bib}

\setlength{\belowdisplayskip}{10pt} \setlength{\belowdisplayshortskip}{10pt}
\setlength{\abovedisplayskip}{10pt} \setlength{\abovedisplayshortskip}{10pt}
\setlength{\intextsep}{20pt}

\end{document}